\documentclass{amsart} 
\usepackage{amsmath}
\usepackage{paralist}
\usepackage{amssymb}
\usepackage{breqn}
\usepackage{titlesec}
\usepackage[misc]{ifsym}
\usepackage{epsfig} 
\usepackage{epstopdf} 
\usepackage[colorlinks=true]{hyperref} 
\hypersetup{urlcolor=blue, citecolor=red}
\allowdisplaybreaks
\titleformat{\section}
  {\normalfont\bfseries\centering}
  {\thesection.}{1em}{}
\usepackage[a4paper, left=1in, right=1in, top=1in, bottom=1in]{geometry}
\theoremstyle{definition}

\newcommand{\ie}{i.e.}
\def\e{\epsilon}

\title[$p-$Laplacian Equations]{Existence of an Infinite Number of Solutions to a Singular Superlinear $p$-Laplacian Equation on Exterior Domains}

\author [ Md Suzan Ahamed, Joseph Iaia ]{}
\subjclass{Primary: 34B40; Secondary: 35B05.}
\keywords{exterior domains, semilinear, superlinear, radial, sign-changing.}

\begin{document}
\maketitle
\centerline{\scshape
Md Suzan Ahamed$^{{\href{mailto:mdsuzanahamed@my.unt.edu}{\textrm{\Letter}}}*1},$ \scshape Joseph Iaia$^{{\href{mailto:iaia@unt.edu}{\textrm{\Letter}}}*1}$}
\medskip
{\footnotesize \centerline{$^1$University of North Texas, USA} } 
\medskip

\begin{abstract} 
\noindent In this paper, we prove the existence of an infinite number of radial solutions of the $p$-$Laplacian$ equation $\Delta_p u + K(|x|) f(u) =0$ on the exterior of the ball of radius $R>0$ in ${\mathbb R}^{N}$ such that $u(|x|)\to 0$ as $|x|\to \infty$ where $f$ grows superlinearly at infinity and is singular at $0$ with $f(u) \sim -\frac{1}{|u|^{m-1}u}$ and $0<m<1$ for small $u$. We also assume $K(|x|) \sim |x|^{-\alpha}$ for large $|x|$ where $N + \frac{m(N-p)}{p-1}< \alpha<2(N-1).$
\end{abstract}

\section{Introduction}
 In this paper, we are interested in radial solutions of
\begin{equation}   \Delta_p u + K(|x|) f(u)=0 \   \textrm{ on } {\mathbb R}^N\backslash B_{R}\label{DE22}  \end{equation}
\begin{equation}  u = 0 \textrm{ on } \partial B_R, \  u\to 0 \textrm{ as } |x| \to \infty \label{DE23} \end{equation}
when $N>2$ where $\Delta_p u=\nabla \cdot (|\nabla u|^{p-2}\nabla u)$ is the $p$-$Laplacian$ with $ 1<p<N,$  and $B_R$ is the ball of radius $R>0$ centered at the origin. \\
Assuming $u(x) =u(|x|) = u(r)$, then $(\ref{DE22})$-$(\ref{DE23})$ becomes
\begin{equation}
|u'|^{p-2}\left((p-1)u'' + \frac{N-1}{r} u'\right) + K(r) f(u) =0  \quad  \textrm{for } R<r < \infty, \label{DE1}
\end{equation}
\begin{equation}   u(R) = 0, \lim_{r \to \infty} u(r) = 0.     \label{DE100}   \end{equation}

\vskip .1 in

We assume $f: {\mathbb R} \backslash\{0\} \to {\mathbb R}$ is odd, locally Lipschitz, and there exists a locally Lipschitz function $g_1:{\mathbb R}\to {\mathbb R}$  such that
\begin{equation}
    f(u)=|u|^{l-1}u + g_1(u)   \textrm{ with } l>p-1 \textrm{ and } 1<p<N \textrm{ for large } |u|  \textrm{ and } \lim\limits_{u \to \infty}\frac{ |g_1(u)|}{|u|^{l}}= 0. \tag{\it H1}
\end{equation}
Thus $f$ is superlinear at $\infty,$ that is,  $\lim\limits_{u\to\infty}\frac{ |f(u)| }{|u|^{p-1}}=\infty$. We also assume there exists a locally Lipschitz function $g_2:{\mathbb R}\to {\mathbb R}$  such that
\begin{equation}
    f(u)=-\frac{1}{|u|^{m-1}u} + g_2(u)   \textrm{ with } 0<m<1 \textrm{ for small } |u|  \textrm{ and }\  g_2(0)=0. \tag{\it H2}
\end{equation}
 In addition, we assume
\begin{equation} f \textrm{  has a unique positive zero, } \beta, \textrm{ with } f<0 \textrm{  on }(0, \beta) \textrm{ and }  f>0 \textrm{  on } (\beta, \infty).\tag{{\it H3}}  \end{equation} 
\vskip .1 in
Let $F(u) = \int_{0}^{u} f(t) \, dt.$ Since $f$ is odd,  $F$ is even. Also, since $0<m<1$ by{\it (H2)}, it follows that $F$ is continuous with $F(0)=0$. Also, it follows from {\it (H1)} and {\it (H3)} that
\begin{align} F \textrm{ has a unique positive zero, }  \gamma, \textrm{ where } 0< \beta< \gamma  \textrm{ with } F<0 \textrm{ on } (0, \gamma) \ \textrm{ and }  F>0 \textrm{  on } (\gamma, \infty).  \label{boss} 
\end{align}
Then, it follows from {\it (H3)} that there exists $F_0>0$ and $f_0>0$ such that
\begin{equation} F(u)\geq -F_0    \textrm{ and } |u|^{m-1}uf(u) \geq -f_0 \textrm{ for all } u \in {\mathbb R\backslash \{0\}}. \label{boss2} \end{equation}
\vskip .1 in
We also assume $K>0$ and $K'$ are continuous on $[R, \infty),$ and there exist constants $K_0, K_1>0$ such that
\begin{align}
\frac{rK'}{K} > -\frac{(N-1)p}{p-1}, \text{ and } \frac{K_0}{r^{\alpha}} \leq  K\leq \frac{K_1}{r^{\alpha_1}}    \textrm{ with }  \ N + \frac{m(N-p)}{p-1}<\alpha_1 \leq \alpha<2(N-1). \tag{\it{H4}} 
\end{align}
\newpage
 In this paper, we prove the following: 
\vskip .1 in

{\bf  THEOREM 1:} {\it  Let $N>2$ and assume (H1)-(H4). Then there exists a non-negative integer $n_0$ such that for each $n\geq n_0$ there is a solution, $u_n\in C^1[R, \infty),$ of (\ref{DE1})-(\ref{DE100}) with exactly $n$ zeros on $(R,\infty)$.} 
 \vskip .1 in
In an earlier paper \cite{ahamed001}, we proved this theorem in the case where $p=2$. In this paper, we extend this result to the $p$-$Laplacian$ with $1<p<N.$ 
\vskip .1 in
A number of papers have been written on this and similar topics. Some have used sub/super solutions or degree theory to prove the existence of a positive solution \cite{castro, chhetri, lee1, lee2, sankar}. Others with more assumptions have been able to prove the existence of an infinite  number of solutions \cite{berestycki1, berestycki2, iaia1, iaia2, iaia3, jones, joshi, mcleod, strauss}.
\vskip .1 in
The main challenge in this paper is the interplay between the singularity of $f(u)$ at $u=0$ and the behavior of $K$ at $r=\infty$. Since it is not immediately clear that solutions of $(\ref{DE22})$-$(\ref{DE23})$ have a finite number of zeros at $r=\infty$ and that the zeros of $u$ are isolated, some care is required in establishing these results. Many papers have been written about $\it{positive}$ solutions, but in this paper, we obtain infinitely many $\it{sign}$-$\it{changing}$ solutions.
\vskip .1 in
We note that at points $r_0$ where $R<r_0<\infty$ and $u(r_0)=0$ then $f(u(r_0))$ is undefined and hence $u$ cannot be twice differentiable at $r_0$. At such a point we will see that $u$ is differentiable and so by a solution of $(\ref{DE1})$-$(\ref{DE100})$ we mean $u_n\in C^1[R, \infty)$ such that $\frac{1}{|u|^m}\in L_{loc}^1[R, \infty)$ and 
$$r^{\frac{N-1}{p-1}} u'=R^{\frac{N-1}{p-1}}u'(R)-\int_{R}^r t^{\frac{N-1}{p-1}} K(t) f(u(t))\, dt  \quad \text{for} \ r\geq R, $$
$$ u(R) = 0, \ \ \lim_{r \to \infty} u(r) = 0.$$
\vskip .1 in

\section {Preliminaries for Theorem 1}

Let $N>2$ and $u$ solve $(\ref{DE1})$-$(\ref{DE100})$ and suppose {\it(H1)}-{\it (H4)} are satisfied. Now let $u(r) = v(r^{\frac{p-N}{p-1}}).$ Then we see (\ref{DE1})-(\ref{DE100}) is equivalent to

\begin{equation}
(p-1)|v'(t)|^{p-2}v''(t) + h(t) f(v(t)) = 0  \   \textrm{ for } 0<t<R^{\frac{p-N}{p-1}},  \label{e9} \end{equation}
\begin{equation} v(0) =0, \  v(R^{\frac{p-N}{p-1}})=0 \label{e10}  \end{equation}

where \begin{equation} t=r^{\frac{p-N}{p-1}} \  \text{and} \   h(t) = \frac{t^{\frac{p(N-1)}{p-N}}K(t ^{\frac{p-1}{p-N}})}{\left(\frac{N-p}{p-1}\right)^p}>0. \label{h eq}\end{equation}
It follows from {\it (H4)} that
\begin{equation}
    h'<0. \label{hprime}
\end{equation}
We now define $\tilde \alpha = \frac{p(N-1)-\alpha (p-1)}{N-p}$ and $\tilde \alpha_1 = \frac{p(N-1)-\alpha_1(p-1)}{N-p}.$ Note that the assumption $N + \frac{m(N-p)}{p-1}<\alpha_1 \leq \alpha<2(N-1)$ from {\it (H4)} implies
\begin{equation}
   \tilde\alpha + m \leq  \tilde \alpha_1 +m<1. \label{alphatilde}
\end{equation}
It also follows from {\it (H4)} and  (\ref{h eq}) that there exist \ $h_0>0 \ \text{and } \  h_1>0$ such that
\begin{equation}    \  h_0 t^{-\tilde \alpha} \leq h \leq h_1 t^{-\tilde \alpha_1}  \textrm{ for }  0< t< R^{\frac{p-N}{p-1}}.  \label{joni} 
\end{equation}
From this and $(\ref{alphatilde})$ it follows that $\frac{h(t) }{t^m}$ and $h(t)$  are integrable on $(0, R^{\frac{p-N}{p-1}}].$
\vskip .1 in
We now consider the following initial value problem
\begin{equation} 
(|v_a'(t)|^{p-2}v_a'(t))' + h(t) f(v_a(t))= 0   \textrm{ for } 0<t<R^{\frac{p-N}{p-1}}, \label{e7} 
\end{equation}
\begin{equation} v_a(0) =0, \ v_a'(0)=a> 0. \label{e8}  \end{equation}
Note: $(|v_a'|^{p-2}v_a')'=(p-1)|v_a'|^{p-2}v_a''.$

\vskip .2 in
{\bf Lemma 2.1:} Let  $N>2$,  $a>0,$  and assume {\it(H1)}-{\it (H4)}. Then there exists an $\e >0$ such that there is a unique solution of  (\ref{e7})-(\ref{e8}) on $[0, \e].$

\vskip .1 in

{\bf Proof:} We first assume $v_a$ solves $(\ref{e7})$-$(\ref{e8})$, then integrate (\ref{e7}) on $(0, t)$ and use (\ref{e8}) along with {\it (H2)} to obtain
\begin{equation} |v_a'(t)|^{p-2}v_a'(t) =a^{p-1} -\int_0^{t} h(x)f(v_a(x))  \, dx  . \label{thomas} \end{equation}
Now define $\Phi_p(x)=|x|^{p-2}x$ for $x\in \mathbb R$, and $p>1.$ Denote its inverse by $\Phi_{p'}(x)$ where $\frac{1}{p}+\frac{1}{p'}=1, \ i.e. \ p'=\frac{p}{p-1}.$ Note that $\Phi_{p}$ and $\Phi_{p'}$ are continuous for $p>1.$ Then from $(\ref{thomas}),$ we obtain
\begin{equation} \Phi_p(v_a'(t)) =a^{p-1} -\int_0^{t} h(x)f(v_a(x))  \, dx  . \label{thomas002} \end{equation}
Therefore,
\begin{equation} v_a'(t) =\Phi_{p'}\left(a^{p-1} -\int_0^{t} h(x)f(v_a(x))  \, dx\right) . \label{thomas003} \end{equation}
Now, integrate (\ref{thomas003})  on $(0,t)$ to obtain
\begin{equation} v_a(t) =\int_0^t\Phi_{p'}\left(a^{p-1} -\int_0^{s} h(x)f(v_a(x))  \, dx\right) \, ds . \label{jeter} \end{equation}
\vskip .1 in
Next, let $w(t)=\frac{v_a(t)}{t}.$ Then, we obtain $w(0)=a>0$ and
\begin{equation} 
w(t) = \frac{a}{t}\int_0^t\Phi_{p'}\left(1 -\frac{1}{a^{p-1}}\int_0^{s} h(x)f(xw(x))  \, dx\right) \, ds .  \label{w eqn}  \end{equation}
Now using the contraction mapping principle, we prove the existence of a unique solution of  (\ref{e7})-(\ref{e8}) on $[0, \epsilon_0]$ for some $\epsilon_0>0.$ Let $\epsilon>0$ and let
$$ \mathcal{B}_\epsilon:=\{ w\in C[0, \epsilon] \ : \ w(0)=a>0 \ \text{and} \  \left|w(t)-a\right|\leq \frac{a}{2} \ \text{on} \  [0, \epsilon]\},$$ where $C[0,\epsilon]$  is the set of all continuous functions on $[0,\epsilon],$ and
$$||w||=\sup_{x\in [0,\epsilon]} |w(x)|.$$ Then $(\mathcal{B}_\epsilon, ||.||)$ is a Banach space \cite{book2010}. Now define a map $T$ on $\mathcal{B}_\epsilon$ by
$$Tw(t)= \begin{cases}
      a    \ \ \quad\quad\quad \quad \quad \quad \quad \  \  \ \ \ \ \ \  \  \quad  \quad\quad\quad  \quad\quad\quad\quad\quad    \quad   \ \   \          \textrm { for } t=0\\
      \frac{a}{t}\int_0^t\Phi_{p'}\left(1 -\frac{1}{a^{p-1}}\int_0^{s} h(x)f(xw(x))  \, dx\right) \, ds   \quad \quad\quad \quad \textrm{\ for } \   0<t \leq \epsilon. 
         \end{cases}$$
Next, note that for any $w\in \mathcal{B}_\epsilon,$ we have
\begin{equation}
    0<\frac{a}{2}\leq w(x)\leq \frac{3a}{2} \quad \text{on} \ [0, \epsilon].\label{wb001}
\end{equation}
Thus
\begin{equation}
    \left|\frac{-1}{x^mw^m(x)}\right|\leq \frac{2^m }{a^m}x^{-m} \quad \text{on} \  (0, \epsilon ]. \label{xw0001}
\end{equation}
Since $g_2$ is locally Lipschitz by $(\it{H2})$ and $g_2(0)=0,$ therefore there exists $L_2>0$ such that
\begin{equation}
    \left|g_2(xw(x))\right|\leq L_2|xw(x)|\leq \frac{3aL_2}{2} |x| \quad \text{on} \  \left [ 0, \epsilon \right ] \ \textrm{ if } \epsilon \textrm{ sufficiently small}. \label{g1}
\end{equation}
Then it follows from $(\ref{joni})$, $(\ref{xw0001})$, $(\ref{g1})$, and $(\it{H2})$ that
$$|h(x)f(xw(x))|\leq \frac{2^mh_1}{a^m}x^{-(\tilde{\alpha}+m)}+\frac{3ah_1L_2}{2}x^{1-\tilde{\alpha}} \quad \text{on} \  \left ( 0, \epsilon \right ] \ \textrm{ if } \epsilon \textrm{ sufficiently small}.$$
Integrating on $(0, t)\subset (0, \epsilon)$, and using $(\ref{alphatilde})$,  we obtain
\begin{equation}
    \int_0^t |h(x)f(xw(x))| \ dx \leq C_1 t^{1-\tilde{\alpha}-m}+C_2 t^{2-\tilde{\alpha}}=t^{1-\tilde{\alpha}-m}[C_1 +C_2 t^{1+m}]\to 0 \quad \text{as} \ t\to 0^+, \label{hf1}
\end{equation}
where $C_1:=\frac{2^mh_1}{a^m(1-\tilde{\alpha}-m)}$ and $C_2:=\frac{3ah_1L_2}{2(1-\tilde{\alpha})}.$ \\
Integrating $(\ref{hf1})$ again, we obtain
\begin{equation}
    \frac{1}{t}\int_0^{t}\int_0^{s} |h(x)f(xw(x))|  \, dx \, ds \to 0 \quad \text{as} \ t\to 0^+.\label{hf3}
\end{equation}
Since $\Phi_{p'}$ is continuous, it follows from $(\ref{hf1})$ that 
\begin{equation}
    \lim_{t\to 0^+}Tw(t)=a=Tw(0) \label{tw0001}
\end{equation}
and so $Tw$ is continuous at $t=0.$  For any $w\in \mathcal{B}_\epsilon$, $Tw(0)=a$ and it is clear from the definition, $(\ref{hf1})$, and $(\ref{tw0001})$ that $Tw(t)$ is continuous for any $t\in [0, \epsilon]$ if $\epsilon$ sufficiently small. 

\vskip .1 in
Note that $\Phi_{p'}'(x)=(p'-1)|x|^{p'-2}=\frac{1}{p-1}|x|^{\frac{2-p}{p-1}}$ and $\Phi_{p'}'(1)=\frac{1}{p-1}.$ Also $|\Phi_{p'}'(x)|\leq \frac{C_p}{p-1}$ if $|x-1|\leq \frac{1}{2},$ where 
$$C_p=\begin{cases}
    \left(\frac{3}{2}\right)^{\frac{2-p}{p-1}}  & \textrm{ when } 1<p\leq 2\\
    \left(\frac{1}{2}\right)^{\frac{2-p}{p-1}}  &  \textrm{ when } p>2.
\end{cases}$$
Thus if $0<t_i<\frac{1}{2} \textrm{ where } i=1, 2$, then by the mean value theorem, there is a $c$ between $1-t_1$ and $1-t_2$ such that
\begin{equation}
    |\Phi_{p'}(1-t_1)-\Phi_{p'}(1-t_2)|=|\Phi_{p'}'(c)||t_1-t_2|\leq \frac{C_p}{p-1} |t_1-t_2|, \quad 0<t_i<\frac{1}{2}. \label{phi001}
\end{equation}
Then $Tw-a=\frac{a}{t}\int_0^t\left[\Phi_{p'}\left(1 -\frac{1}{a^{p-1}}\int_0^{s} h(x)f(xw(x))  \, dx\right)-1\right] \, ds$, and thus
$$|Tw-a |\leq \frac{a}{t}\int_0^t\left|\Phi_{p'}\left(1 -x_*\right) \,-1\right|\, ds \quad \textrm{where } x_*=\frac{1}{a^{p-1}}\int_0^{s} h(x)f(xw(x))  \, dx.$$
Now using $(\ref{phi001})$ in the above inequality, we obtain
$$|Tw-a | \leq \frac{C_p}{(p-1)a^{p-2}t}\int_0^t\int_0^{s} |h(x)f(xw(x))|  \, dx\, ds.$$
Then it follows from $(\ref{hf3})$ that
$$|Tw-a | \leq\frac{C_p}{(p-1)a^{p-2}t}\int_0^t\int_0^{s} |h(x)f(xw(x))|  \, dx\, ds \to 0 \quad \text{as} \ t\to 0^+.$$
Thus $|Tw(t)-a|\leq \frac{a}{2}$ for sufficiently small $t>0.$ Therefore $T:\mathcal{B}_\epsilon\longrightarrow \mathcal{B}_\epsilon\text{ if } \epsilon$ is sufficiently small.
\vskip .1 in
We next prove that $T$  is a contraction mapping if $\epsilon$ is sufficiently small. \\

Let $w_1,w_2 \in \mathcal{B}_\epsilon$. Then
$$    |Tw_1- Tw_2|\leq \frac{a}{t}\int_0^{t}\left|\Phi_{p'}\left(1 -\frac{1}{a^{p-1}}\int_0^{s} h(x)f(xw_1(x))  \, dx\right)-\Phi_{p'}\left(1 -\frac{1}{a^{p-1}}\int_0^{s} h(x)f(xw_2(x))  \, dx\right)\right| \, ds.$$
It follows from $(\ref{phi001})$ that
\begin{equation}
    |Tw_1- Tw_2|\leq\frac{C_p}{(p-1)a^{p-2}t}\int_0^{t}\int_0^{s} |h(x)[f(xw_1(x))-f(xw_2(x)]|  \, dx \, ds.  \label{tw002}
\end{equation}
Note that by {\it (H2)}
$$| f(xw_1)-f(xw_2)|\leq \frac{1}{x^m} \left|\frac{1}{{w_1}^m}-\frac{1}{{w_2}^m}\right|+L_2  x |w_1-w_2|,
$$
where $L_2$ is the Lipschitz constant for $g_2$ near $x=0.$ Then using the mean value theorem, for any $x\in (0, \epsilon)$ there is a $c_x$ with $\frac{a}{2}\leq c_x \leq \frac{3a}{2}$ such that
$$ \frac{1}{x^m}\left|\frac{1}{{w_1}^m}-\frac{1}{{w_2}^m}\right|=\frac{m}{x^m c_x^{m+1}}|w_1-w_2|\leq \frac{2^{m+1}m}{a^{m+1}} |w_1-w_2|x^{-m} \leq \frac{2^{m+1}m}{a^{m+1}} ||w_1-w_2||x^{-m}.$$
It follows from this, $(\ref{alphatilde}), (\ref{joni})$, and integrating twice on $(0, t)$ that
\[\int_0^t \int_0^s|h(x)(f(xw_1)-f(xw_2))| \ dx \ ds\leq \left[\frac{2^{m+1}m h_1 t^{2-(\tilde{\alpha_1}+m)}}{(1-\tilde{\alpha_1}-m)(2-\tilde{\alpha_1}-m)a^{m+1}} + \frac{L_2 h_1 t^{3-\tilde{\alpha_1}}}{(2-\tilde{\alpha_1})(3-\tilde{\alpha_1})}\right]||w_1-w_2||.\]
Thus
\begin{equation}
    \frac{2}{(p-1)a^{p-2}t}\int_0^t \int_0^s|h(x)(f(xw_1)-f(xw_2))| \ dx \ ds\leq (C_3 \epsilon^{1-(\tilde{\alpha_1}+m)}+ C_4 \epsilon^{2-\tilde{\alpha_1}})||w_1-w_2|| \quad \text{on} \ (0, \epsilon), \label{nr1}
\end{equation}
where $C_3=\frac{2^{m+2}m h_1}{(p-1)(1-\tilde{\alpha_1}-m)(2-\tilde{\alpha_1}-m)a^{p+m-1}}$ and $C_4=\frac{2L_2 h_1}{(p-1)(2-\tilde{\alpha_1})(3-\tilde{\alpha_1})a^{p-2}}.$\\
Using $(\ref{nr1})$ in $(\ref{tw002})$, we see that
$$|Tw_1- Tw_2|\leq C_\epsilon ||w_1-w_2||, \quad \textrm{where } C_\epsilon :=C_3  \epsilon^{1-(\tilde{\alpha_1}+m)}+ C_4 \epsilon^{2-\tilde{\alpha_1}}.$$
Note that $C_\epsilon\to 0$ as $\epsilon\to 0^+.$ Now choose $\epsilon_0$ such that $C_{\epsilon_0}<1$, then $T$ is a contraction on $\mathcal{B}_{\epsilon_0}.$ Thus, by the contraction mapping principle \cite{book2010}, there exists a unique $w\in \mathcal{B}_{\epsilon_0}$ such that $Tw=w$ on $[0, \epsilon_0]$, and hence $v_a(t)=tw(t)$ is a solution of $(\ref{e7})$-$(\ref{e8})$ on $[0, \epsilon_0].$ \\

Moreover, from $(\ref{wb001}),$ we see that $\frac{a}{2}t\leq v_a\leq \frac{3a}{2}t.$ So $v_a$ is bounded on $[0, \epsilon_0].$ Also, from $(\ref{thomas003})$ and $(\ref{hf1})$ and since $\Phi_{p'}$ is continuous, it follows that $v_a'$ is continuous as well as bounded on $[0, \epsilon_0].$ This completes Lemma 2.1. \qed

\vskip .2 in

{\bf Lemma 2.2:} Let $N>2, \ a>0$, and assume {\it (H1)}-{\it (H4)}. Suppose $v_a$ solves (\ref{e7})-(\ref{e8})  on $[0, z_a]$ where $0<z_a<R^{\frac{p-N}{p-1}}.$ If $v_a(z_a)=0,$ then $v_a'(z_a)\neq 0.$ 
\vskip .1 in

{\bf Proof:} Let
$$E(t):=\frac{p-1}{p}\frac{|v_a'(t)|^p}{h(t)}+F(v_a(t)).$$
Since $h$ is decreasing (from $(\ref{hprime})$), it follows from $(\ref{e7})$ that 

\begin{equation}
    E'(t)=-\frac{p-1}{p}\frac{|v_a'(t)|^ph'(t)}{h(t)^2}\geq 0 \quad \textrm{for all } 0<t\leq z_a. \label{energy}
\end{equation} 
Note that $E(0)=0$. It follows from $(\ref{energy})$ that if $v_a(z_a)= v_a'(z_a)=0$ with $0<z_a<R^{\frac{p-N}{p-1}}$, then $E(z_a)=0$. But since $E$ is non-decreasing, this implies $v_a\equiv v_a'\equiv0$ on $[0, z_a]$ which contradicts that $v_a'(0)=a>0.$ Thus if  $v_a(z_a)=0,$ then $v_a'(z_a)\neq 0.$ \qed

\vskip .2 in

{\bf Lemma 2.3:}  Let $N>2, a>0$, and assume {\it (H1)}-{\it (H4)}. Let $v_a$ solve (\ref{e7})-(\ref{e8}), then $v_a$ and $v_a'$ are defined and continuous on $[0, R^{\frac{p-N}{p-1}}]$, and they vary continuously with $a$.  In addition, $v_a$ has at most a finite number of zeros on $[0, R^{\frac{p-N}{p-1}}]$.

\vskip .1 in

{\bf Proof:} Let $t_0>0$. We first show that solutions of $(\ref{e7})$ with $v_a(t_0)=b, v_a'(t_0)=c$ where $|b|+|c|>0$ are unique in a neighborhood of $t_0.$

Suppose $v_{1, a}$ and $v_{2, a}$ solve $(\ref{e7})$ and 
$$v_{1, a}(t_0)=v_{2, a}(t_0)=b,$$
$$v_{1, a}'(t_0)=v_{2, a}'(t_0)=c.$$
\textbf{Case I:} $b\neq 0$
\vskip .1 in
Integrating $(\ref{e7})$ on $(t_0, t)$ gives
$$v_{1,a}'(t) =\Phi_{p'}\left(b^{p-1} -\int_{t_0}^{t} h(x)f(v_{1, a}(x))  \, dx\right),$$
$$v_{2,a}'(t) =\Phi_{p'}\left(b^{p-1} -\int_{t_0}^{t} h(x)f(v_{2, a}(x))  \, dx\right).$$
Since $\Phi_{p'}'(x)$ is continuous at $b\neq 0$, then there exists $C_5$ such that $|\Phi_{p'}'(x)|\leq C_5$ in a neighborhood of $b.$ Then by the mean value theorem, there exists a $C_a$ with $C_a$ close to $b\neq 0$ such that
$$|(v_{1,a}(t)-v_{2,a}(t))'|\leq |\Phi_{p'}'(C_a)|\int_{t_0}^{t} h(x)|f(v_{1, a}(x))-f(v_{2, a}(x))|  \, dx\leq C_5 \int_{t_0}^{t} h(x)|f(v_{1, a}(x))-f(v_{2, a}(x))|  \, dx$$
Similarly, since $f$ is continuous at $b\neq 0$, then there exists $C_f$ such that $|f(x_1)-f(x_2)|\leq C_f|x_1-x_2|$ in a neighborhood of $b.$ Then we see that
$$|(v_{1,a}(t)-v_{2,a}(t))'|\leq C_5C_f\int_{t_0}^{t} h(x)|v_{1, a}(x)-v_{2, a}(x)|  \, dx.$$
Integrating this on $(t_0, t)$, we obtain
\begin{equation}
   | v_{1,a}(t)-v_{2,a}(t)|\leq \int_{t_0}^{t}|(v_{1,a}(x)-v_{2,a}(x))'| \, dx\leq C_5C_f\int_{t_0}^{t} \int_{t_0}^{s} h(x)|v_{1, a}(x)-v_{2, a}(x)|  \, dx \, ds. \label{eq001}
\end{equation}
Notice
\begin{equation}
\int_{t_0}^{t} \int_{t_0}^{s} h(x)|v_{1, a}(x)-v_{2, a}(x)|  \, dx \, ds\leq \int_{t_0}^{t} \int_{t_0}^{t} h(x)|v_{1, a}(x)-v_{2, a}(x)|  \, dx \, ds\leq R^{\frac{p-N}{p-1}}\int_{t_0}^{t} h(x)|v_{1, a}(x)-v_{2, a}(x)|  \, dx.\label{aaa001}
\end{equation}
It follows from $(\ref{eq001})$-$(\ref{aaa001})$ that
\begin{equation}
    |v_{1,a}(t)-v_{2,a}(t)|\leq R^{\frac{p-N}{p-1}}C_5C_f\int_{t_0}^{t} h(x)|v_{1, a}(x)-v_{2, a}(x)|  \, dx.\label{eq003}
\end{equation}
Then $v_{1, a}-v_{2, a}\equiv0$ in a neighborhood of $t_0$ by the Gronwall Inequality \cite{book1969}, and hence $v_{1, a}\equiv v_{2, a}$ in a neighborhood of $t_0.$
\vskip .1 in
\textbf{Case II:} $b=0$
\vskip .1 in
It follows from Lemma 2.2 that $v_{1, a}'(t_0)= v_{2, a}'(t_0)\neq 0.$ In this case, we can apply the contraction mapping principle to get a $\it{unique}$ solution of $(\ref{e7})$, $v_*,$ with $v_*(t_0)=0, v_*'(t_0)=v_{1, a}'(t_0)= v_{2, a}'(t_0)$. Thus $ v_{1, a}\equiv v_{2, a}\equiv v_*$ on $(t_0-\epsilon, t_0+\epsilon)$, and hence $v_a$ is $\it{unique}$ in a neighborhood of $t_0$ where $0<t_0<R^{\frac{p-N}{p-1}}.$ 

\vskip .1 in
Let $[0, d)$ be the maximal half-open interval of existence of the solution of $(\ref{e7})$-$(\ref{e8})$ . Now, note it follows from $(\ref{e7})$-$(\ref{e8})$ that
\begin{equation}
    \left(\frac{p-1}{p}|v_a'(t)|^p+h(t)F(v_a(t))\right)'= h'(t)F(v_a(t)) \label{en001}.
\end{equation}
Integrating both sides on $[\epsilon_0, t]$ with $t\leq d$ where $[0, \epsilon_0]$ is the interval of existence obtained in Lemma 2.1, we obtain
$$\frac{p-1}{p}|v_a'(t)|^p+h(t)F(v_a(t))= \frac{p-1}{p}|v_a'(\epsilon_0)|^p+h(\epsilon_0)F(v_a(\epsilon_0))+\int_{\epsilon_0}^t h'(s)F(v_a(s)) \ ds.$$
It follows from $(\ref{boss2})$ that $F\geq-F_0$ for some $F_0.$ Since $h>0$, therefore $hF\geq-hF_0, \ h'F\leq-F_0h' \ (\because h'<0 \ \text{from (\ref{hprime})})$, and hence
$$|v_a'(t)|^p\leq \frac{p}{p-1}\left[\frac{p-1}{p}|v_a'(\epsilon_0)|^p+h(\epsilon_0)F(v_a(\epsilon_0))+F_0h(\epsilon_0)\right] \quad \text{on} \ [\epsilon_0, d].$$
We then obtain $|v_a'(t)|\leq C_6 \ \text{on} \ [\epsilon_0, d]$, where
$$ C_6:=\left[|v_a'(\epsilon_0)|^p+\frac{p}{p-1}h(\epsilon_0)(F(v_a(\epsilon_0))+F_0)\right]^{\frac{1}{p}}.$$
Also, since $v_a(t)=v_a(\epsilon_0)+\int_{\epsilon_0}^t v_a'(s) \ ds$, then $|v_a(t)|\leq C_7 \ \text{on} \ [\epsilon_0, d]$ where $C_7:=C_5(R^{\frac{p-N}{p-1}}-\epsilon_0)+|v_a(\epsilon_0)|.$
\vskip .1 in
Now combining this with the bounds obtained for $v_a$ and $v_a'$ on $[0, \epsilon_0]$ in Lemma 2.1, we see that $v_a$ and $v_a'$ are uniformly bounded on $[0,d]$ by a constant that only depends on $\epsilon_0$ and $R.$ It follows from this and the earlier proof of uniqueness that $v_a$ and $v_a'$ can be uniquely extended to be defined on $[0, R^{\frac{p-N}{p-1}}]$,  and $v_a$, $v_a'$ are uniformly bounded on $[0, R^{\frac{p-N}{p-1}}]$. It also follows from this that $v_a$ varies continuously with $a>0$.
\vskip .1 in

We now prove that $v_a$ has only a finite number of zeros on $[0, R^{\frac{p-N}{p-1}}]$. Suppose by way of contradiction that $v_a$ has an infinite number of zeros, $z_k$, with $0<z_k<z_{k+1}.$ Then $v_a$ would have an infinite number of local extrema, $M_k$, with $z_k < M_{k+1}<z_{k+1}$. Since $z_k$ and $M_k$ are in the compact set $[0, R^{\frac{p-N}{p-1}}],$ it follows that $z_k$ would have a  limit point, $z^*$, and $v_a(z^*)=0 \ (\because v_a(z_k)\to v_a(z^*) \text{ and } v_a(z_k)=0 \ \forall k)$. Also, we would have $M_k\to z^*,$  and $v_a'(z^*)=0 \ (\because v_a'(M_k)\to v_a'(z^*) \text{ and } v_a'(M_k)=0 \ \forall k)$ which contradicts Lemma 2.2. Thus, $v_a$ has at most a finite number of zeros on $[0, R^{\frac{p-N}{p-1}}]$.\qed

\vskip .2 in

{\bf Lemma 2.4:} Let $N>2,$ $a>0$, and assume {\it (H1)}-{\it (H4)}. Suppose $v_a$ solves (\ref{e7})-(\ref{e8}). Then
$$\max\limits_{\left[0, \ R^{\frac{p-N}{p-1}}\right]} v_a \to \infty \  \text{as} \ a \to \infty.$$
\vskip .1 in

{\bf Proof:} First, we prove the statement for the interval $[0, t_0]$ where $t_0>0$ is sufficiently small and $t_0<R^{\frac{p-N}{p-1}}$. Suppose by way of contradiction that there is a $C_8>0$ such that $|v_a(t)|\leq C_8$ for all $t\in [0, \ t_0]$ and for all $a>0.$ Then there exists $C_9>0$ such that $|g_2(v_a)|\leq C_9$ for all $t\in [0, \ t_0]$.  
\vskip .1 in
Integrating $(\ref{e7})$ on $(0, t)\subset (0, t_0),$ and using {\it (H1)} and $(\ref{joni})$, we get
$$v_a'(t)^{p-1} \geq a^{p-1} -\int_0^{t} h(x)\left[\frac{-1}{v_a^m} +g_2(v_a)\right]  \, dx \geq a^{p-1}-\int_0^{t} |h(x)g_2(v_a)|  \, dx \geq a^{p-1}-h_1C_9\int_0^{t}  x^{-\tilde \alpha_1} \, dx.$$
Thus
$$v_a'(t)^{p-1} \geq  a^{p-1}-h_1C_9\frac{R^{\frac{p-N}{p-1}(1-\tilde \alpha_1)}}{1-\tilde \alpha_1}= a^{p-1}-C_{10},\quad \textrm{where  } C_{10}:=h_1C_9\frac{R^{\frac{p-N}{p-1}(1-\tilde \alpha_1)}}{1-\tilde \alpha_1}.$$

It follows from this that if $a>0$ is sufficiently large then $a^{p-1}-C_{10}>0$ and 
\begin{equation}
    v_a'(t) \geq (a^{p-1}-C_{10})^{\frac{1}{p-1}}\to \infty \textrm{ as } a\to \infty. \label{va001}
\end{equation}
Now integrating $(\ref{va001})$ on $[0, t]\subset [0, t_0],$ we get
\begin{equation}
    v_a(t)=\int_0^{t} v_a'(s) \, ds \geq \int_0^t (a^{p-1}-C_{10})^{\frac{1}{p-1}} \, ds=(a^{p-1}-C_{10})^{\frac{1}{p-1}}t. \label{abc001}
\end{equation}
Evaluating this inequality at $t=t_0$ and using the assumption $v_a\leq C_8$, we see that
$$C_8\geq v_a(t_0)\geq (a^{p-1}-C_{10})^{\frac{1}{p-1}}t_0\to \infty \quad \textrm{as } a\to \infty.$$
This is a contradiction, and hence
$$\max\limits_{\left[0, \ t_0\right]} v_a \to \infty \  \text{as} \ a \to \infty.$$
Now, since $\max\limits_{t\in [0, \ R^{\frac{p-N}{p-1}}]} v_a(t)\geq \max\limits_{t\in [0, \ t_0]} v_a(t),$ therefore we see 
$$\max\limits_{[0, \ R^{\frac{p-N}{p-1}}]} v_a \to \infty \  \text{as} \ a \to \infty.$$
Therefore, the lemma holds. \qed

\vskip .2 in

{\bf Lemma 2.5:}   Let $N>2, a>0$, and assume {\it (H1)}-{\it(H4)}. Suppose $v_a$ solves (\ref{e7})-(\ref{e8}). Then $v_a$ has a first local maximum, $M_a$,  if $a$ is sufficiently large. In addition,  $M_a \to 0 $  and $v_a(M_a)\to \infty$ as $a \to \infty.$
\vskip .1 in

{\bf Proof:}  Suppose not.  Therefore, assume that $v_a>0 \textrm{ and } v_a'>0$ on $(0, R^{\frac{p-N}{p-1}}]$ for all $a>0.$ 
\vskip .1 in
We first show that for sufficiently large $a>0,$ there exists a $t_{a, \beta}\in (0, R^{\frac{p-N}{p-1}})$ such that $v_a(t_{a, \beta})=\beta$ and $0<v_a<\beta$ on $(0, t_{a, \beta})$. Also we show that $t_{a, \beta}\to 0$ as $a\to \infty$.
\vskip .1 in
Suppose there is no such $t_{a, \beta}$, so assume $v_a<\beta$ on $(0, R^{\frac{p-N}{p-1}}]$ for all $a>0.$  Then $f(v_a)<0$ on $(0, R^{\frac{p-N}{p-1}}]$, and hence from $(\ref{e7})$ we see $v_a'>a \textrm{ for all } t\in (0, R^{\frac{p-N}{p-1}}].$ Integrating this inequality on $(0, t),$ we get
\begin{equation}
    v_a(t)>at \quad \forall   t\in (0, R^{\frac{p-N}{p-1}}] \textrm{ and } \forall  a>0. \label{help001}
\end{equation}
In particular, at $t=\frac{1}{2}R^{\frac{p-N}{p-1}},$ we see that
$$ \beta>v_a(R^{\frac{p-N}{p-1}})>\frac{1}{2}R^{\frac{p-N}{p-1}}a\to \infty \textrm{ as } a\to \infty, \ \  \textrm{a contradiction to our assumption that } v_a<\beta.$$
Therefore, there exists $t_{a, \beta}>0$ such that $v_a(t_{a, \beta})=\beta$ and $0<v_a<\beta$ on $(0, t_{a, \beta})$. Evaluating $(\ref{help001})$ at $t=t_{a, \beta},$ we get
\begin{equation}
    t_{a, \beta}< \frac{\beta}{a}\to 0 \quad \textrm{as } a\to \infty . \label{eee001}
\end{equation}
Next, note it follows from $(\ref{e7})$ that
\begin{equation}
    \left(\left(\frac{v_a'}{v_a}\right)^{p-1}\right)'=-\frac{h(t)f(v_a)}{v_a^{p-1}}-(p-1)\left(\frac{v_a'}{v_a}\right)^p  \quad \textrm{on } [t_{a, \beta}, R^{\frac{p-N}{p-1}}]. \label{vprimev}
\end{equation}
Since we are assuming that $v_a>0 \textrm{ and } v_a'>0$ on $(0, R^{\frac{p-N}{p-1}}]$, then $(\ref{vprimev})$ implies that
$$\left(\left(\frac{v_a'}{v_a}\right)^{p-1}\right)'\leq -(p-1)\left(\frac{v_a'}{v_a}\right)^p \quad \textrm{on } [t_{a, \beta}, R^{\frac{p-N}{p-1}}] .$$
Now let $y(t)=\left(\frac{v_a'(t)}{v_a(t)}\right)^{p-1}.$ Then $y^{\frac{p}{p-1}}(t)=\left(\frac{v_a'(t)}{v_a(t)}\right)^{p},$  and from the above inequality, we see
$$\frac{y'}{y^{\frac{p}{p-1}}}\leq -(p-1) \quad \textrm{on } [t_{a, \beta}, R^{\frac{p-N}{p-1}}].$$
Integrating this on $(t_{a, \beta}, t)$, we obtain $y^{\frac{1}{p-1}}(t)\leq \frac{1}{t-t_{a, \beta}},$ and hence
\begin{equation}
    \frac{v_a'(t)}{v_a(t)}\leq \frac{1}{t-t_{a, \beta}} \quad \textrm{ for all } t\in (t_{a, \beta}, R^{\frac{p-N}{p-1}}]. \label{bound001}
\end{equation}
Again, from $(\ref{vprimev})$ we have
\begin{equation}
    \left(\left(\frac{v_a'}{v_a}\right)^{p-1}\right)'\leq -\frac{h(t)f(v_a)}{v_a^{p-1}} \quad \textrm{on } (t_{a, \beta}, R^{\frac{p-N}{p-1}}]. \label{vprimehf}
\end{equation}
Since we are assuming that $v_a$ is positive and increasing on $(0, R^{\frac{p-N}{p-1}}]$, mimicking the proof of Lemma 2.4, we see that if $t_0>0$ and $0<t_{a, \beta}<t_0<R^{\frac{p-N}{p-1}},$ then $v_a(t_0)\to\infty$ as $a\to \infty$. \\
Now integrating $(\ref{vprimehf})$ on $(t_0, R^{\frac{p-N}{p-1}}],$ and using our assumption $l-p+1>0$ from {\it (H1)}  and $(\ref{bound001}),$ we obtain
$$\left(\frac{v_a'}{v_a}\right)^{p-1}\left(R^{\frac{p-N}{p-1}}\right)\leq\left(\frac{v_a'}{v_a}\right)^{p-1}\left(t_0\right) -h\left(R^{\frac{p-N}{p-1}}\right)v_a^{l-p+1}\left(t_0\right)\left(R^{\frac{p-N}{p-1}}-t_0\right).$$
Since $t_{a, \beta}\to 0$ as $a\to 0$ (from $(\ref{eee001})$), using $(\ref{bound001})$ in the above inequality and Lemma 2.4, we see that
\[ \left(\frac{v_a'}{v_a}\right)^{p-1}\left(R^{\frac{p-N}{p-1}}\right) \leq\left(\frac{1}{t_0-t_{a, \beta}}\right)^{p-1} -h\left(R^{\frac{p-N}{p-1}}\right)v_a^{l-p+1}\left(t_0\right)\left(R^{\frac{p-N}{p-1}}-t_0\right)\to -\infty \quad \textrm{ as } a\to \infty.\]
This is a contradiction to our assumption that both $v_a \textrm{ and } v_a'$ are positive on $(0, R^{\frac{p-N}{p-1}}].$ Therefore $v_a$ must have a local maximum, $M_a$, on $(0, R^{\frac{p-N}{p-1}}]$ if $a$ is sufficiently large. It then follows from Lemma 2.4 that $v_a(M_a)\to \infty$ as $a \to \infty.$

\vskip .1 in
Next, we show $M_a \to 0$ as $a \to\infty$.

\vskip .1 in
Suppose $M_a$ does not go to 0 as $a \to \infty$. Then, there is a $C_M>0$ such that $M_a > C_M$ for large $a$.
Again, since $t_{a,\beta}\to 0$ as $a\to \infty$ (by $(\ref{eee001})$), it follows that $t_{a,\beta}<\frac{M_{a}}{2}< M_a$ for large $a$ and from $(\ref{e7})$ we see that $h(t) f(v_a)>0,$ so $(p-1)|v_a'|^{p-2}v_a''=(|v_a'|^{p-2}v_a')'<0$ and so $v_a''<0$ on $(t_{a,\beta}, M_a)$. Therefore,  using the concavity of $v_a$, we obtain
\begin{equation} v_a\left(\lambda_a t_{a, \beta} + (1-\lambda_a) M_a\right)  \geq  \lambda_a\beta + (1-\lambda_a)v_a(M_a) \quad  \text{for }0<\lambda_a<1.  \label{ringo} \end{equation}
Now choose $\lambda_a$ so that $\lambda_a t_{a, \beta} + (1-\lambda_a) M_a = \frac{M_a}{2}$ $\left( \ie  \ \lambda_a = \frac{M_a}{2(M_a-t_{a,\beta})}=  \frac{1}{2\left(1- \frac{t_{a,\beta}}{M_a}\right)}\right) $, then we see from $(\ref{eee001})$ that $\frac{t_{a,\beta}}{M_a}\leq \frac{t_{a,\beta}}{C_M} \to 0$ and so $\lambda_a \to \frac{1}{2}$  as $a \to \infty$. Then by (\ref{ringo}) and Lemma 2.4,
\begin{align} v_a\left(\frac{M_a}{2}\right) \geq  \lambda_a \beta + (1 - \lambda_a)v_a(M_a) =   \frac{\beta}{2(1-\frac{t_{a,\beta}}{M_a})}  + \left(1-\frac{1}{2(1-\frac{t_{a,\beta}}{M_a})}\right)v_a(M_a) \to \infty \ 
\textrm{ as } a \to \infty.  \label{tutti} 
\end{align}
Next, we  integrate (\ref{e7}) on $(t, M_a)$ with $\frac{M_a}{2} < t <M_a.$ Recalling that $v_a$ is increasing on $(\frac{M_a}{2}, M_a)$, and $v_a\to\infty$  on $(\frac{M_a}{2}, M_a)$ from (\ref{tutti}), it follows from  ({\it H1}) that $ f(v_a) \geq C_{2, f} v_a^l$ for sufficiently large $a$ with $C_{2, f}>0$ on $(\frac{M_a}{2}, M_a)$, and so integrating $(\ref{e7})$ on $[t, M_a]\subset [\frac{M_a}{2}, M_a]$, we obtain
$$  (v_a'(t))^{p-1}= \int_{t}^{M_a}     h(s) f(v_a)  \, ds \geq C_{2, f}\int_{t}^{M_a}  h(s) v_a^l(s)   \, ds \geq  C_{2, f} v_a^{l}(t)\int_{t}^{M_a} h(s)  \, ds   \textrm{  \ \ on  } \left(\frac{M_a}{2}, M_a\right). $$
Thus, $$  \frac{ v_a'}{v_{a}^{\frac{l}{p-1}}} \geq C_{2, f}^{\frac{1}{p-1}}\left(\int_{t}^{M_a} h(s) \, ds\right)^{\frac{1}{p-1}} \textrm{ \ \ on } \left(\frac{M_a}{2}, M_a\right).   $$
Integrating this on $ (t,M_a)$ and recalling from {\it (H1)}, we obtain for any $t\in \left(\frac{M_a}{2}, M_a\right)$
\begin{equation}  v_{a}^{\frac{p-l-1}{p-1}}(t)  \geq v_a^{\frac{p-l-1}{p-1}}(t) - v_a^{\frac{p-l-1}{p-1}}(M_a) \geq \frac{l-p+1}{p-1}C_{2, f}^{\frac{1}{p-1}} \int_{t}^{M_a} \left(\int_{s}^{M_a} h(x)  \, dx \right)^{\frac{1}{p-1}} ds. \label{1111} \end{equation}
Now, we evaluate (\ref{1111}) at $t=\frac{M_a}{2}.$ Since $v_a\left(\frac{M_a}{2}\right)\to \infty$ as $a \to \infty$ from (\ref{tutti}), then
$$  0 \leftarrow   v_a^{\frac{p-l-1}{p-1}}\left(\frac{M_a}{2}\right) \geq \frac{l-p+1}{p-1}C_{2, f}^{\frac{1}{p-1}} \int_{\frac{M_a}{2}}^{M_a} \left(\int_{s}^{M_a} h(x)  \, dx \right)^{\frac{1}{p-1}} ds \quad \text{as} \ a\to \infty.  $$
However,  as $a \to \infty$ the right-hand side goes to a positive  constant  (since $M_a \geq C_M>0$), and the left-hand side goes to 0, so we obtain a contradiction. Thus, $M_a \to 0$ as $a \to \infty$. This completes the proof.\qed

\vskip .2 in

{\bf Lemma 2.6:} Let $N>2, a>0$, and assume {\it (H1)}-{\it (H4)}. Suppose $v_a$ solves (\ref{e7})-(\ref{e8}). Then $v_a$ has a first zero, $z_a$, with $0< z_a <R^{\frac{p-N}{p-1}}$ if $a$ is sufficiently large. In addition, $z_a \to 0$ as $a \to \infty$.

\vskip .1 in

{\bf Proof:} From Lemma 2.5, we know that $v_a$ has a first local maximum, $M_a$, with $0<M_a<R^{\frac{p-N}{p-1}}$ for sufficiently large $a>0.$ Now suppose by way of contradiction that $v_a>0$ on $(M_a, R^{\frac{p-N}{p-1}})$ for all $a$ sufficiently large. From $(\ref{energy})$, we know that $E$ is non-decreasing on  $(0, R^{\frac{p-N}{p-1}}]$. It follows then that
\begin{equation} \frac{p-1}{p}  \frac{|v_{a}'|^p}{h} + F(v_a)  \geq F(v_a(M_a)) \textrm{ \ on } (M_a, R^{\frac{p-N}{p-1}}). \label{frutti} \end{equation}
$\textbf{Claim:}$ For sufficiently large $a>0,$ if $v_a>0$ on $(M_a, R^{\frac{p-N}{p-1}})$ then $v_a'<0$ on $(M_a, R^{\frac{p-N}{p-1}}).$
\vskip .1 in

$\textbf{Proof of Claim:}$ Suppose there is a smallest critical point, $t_a$, with $M_a<t_a<R^{\frac{p-N}{p-1}},$ $v_a'(t_a)=0$ and $v_a(t_a)>0.$ Then from $(\ref{frutti}),$ we have $F(v_a(t_a))  \geq F(v_a(M_a))$. Since $v_a(M_a)\to \infty$ as $a\to \infty$ from Lemma 2.5, therefore $v_a(M_a)>\gamma$ for sufficiently large $a>0.$ So $F(v_a(M_a))>0$, and hence $F(v_a(t_a))>0$ for sufficiently large $a>0$ and thus $v_a(t_a)>\gamma.$ Since $v_a$ decreasing on $(M_a, t_a)$ then $v_a(t)>\gamma>\beta$ on $(M_a, t_a)$, therefore $f(v_a(t))>0$ on $(M_a, t_a),$ and hence $\int_{M_a}^{t_a}hf(v_a) \, ds>0.$ However, if we integrate $(\ref{e7})$ on $(M_a, t_a)$, we get $\int_{M_a}^{t_a}hf(v_a) \, ds=0,$ a contradiction. So $v_a'<0$ on $(M_a, R^{\frac{p-N}{p-1}}).$ This proves the claim.

\vskip .1 in
Next, rewriting $(\ref{frutti})$ gives
\begin{equation}
     \sqrt[p]{h}\leq \sqrt[p]{\frac{p-1}{p}}\frac{| v_a'|}{ \sqrt[p]{F(v_a(M_a)) - F(v_a)}},  \label{rewrite001}
\end{equation}
and integrating on $(M_a, R^{\frac{p-N}{p-1}})$ gives
\begin{equation}  \int_{M_a}^{R^{\frac{p-N}{p-1}}} \sqrt[p]{h(t)} \, dt \leq \sqrt[p]{\frac{p-1}{p}} \int_{v_a\left(R^{\frac{p-N}{p-1}}\right)}^{v_a(M_a)}  \frac{  dt}{ \sqrt[p]{F(v_a(M_a)) - F(t)}}  
\leq \sqrt[p]{\frac{p-1}{p}}  \int_{0}^{v_a(M_a)}  \frac{  dt}{ \sqrt[p]{F(v_a(M_a)) - F(t)}}. \label{diddley}  \end{equation}
Now, letting $t= v_a(M_a)s$ and changing variables on the right-most integral in (\ref{diddley}), we obtain
\begin{equation}
    \sqrt[p]{\frac{p-1}{p}}  \int_{0}^{v_a(M_a)}  \frac{  dt}{ \sqrt[p]{F(v_a(M_a)) - F(t)}}  = \sqrt[p]{\frac{p-1}{p}}\frac{v_a(M_a)}{\sqrt[p]{F(v_a(M_a))}} \int_{0}^1 \frac{ ds }{\sqrt[p]{1 - \frac{F(v_a(M_a)s)}{F(v_a(M_a))} }}.\label{abc00001}
\end{equation}
We know from earlier that $v_a(M_a) \to \infty$ as $a \to \infty$ and so it follows from the superlinearity of $f$  (from {\it (H1)}) that  $\sqrt[p]{\frac{p-1}{p}}\frac{v_a(M_a)}{\sqrt[p]{F(v_a(M_a))}}   \to 0$ as $a\to \infty$. In addition, since $\frac{F(v_a(M_a)s)}{F(v_a(M_a))} \to s^{l+1}$ uniformly on $[0, 1]$ as $a\to \infty$, we obtain
$$\int_{0}^1 \frac{ ds }{    \sqrt[p]{1 - \frac{F(v_a(M_a)s)}{F(v_a(M_a))}       }         }    
\to \int_{0}^{1} \frac{ ds }{\sqrt[p]{1 - s^{l+1} }    } < \infty \textrm{ as }a \to \infty.$$
Thus, it follows from these facts, $(\ref{diddley})$, and $(\ref{abc00001})$ that
\begin{equation}
      \int_{M_a}^{R^{\frac{p-N}{p-1}}} \sqrt[p]{h} \, dt \to 0  \textrm{ as }  a \to \infty \label{sqth1}.
\end{equation}
However, we know $M_a \to 0$ from Lemma 2.5, and so the left-hand side of  (\ref{sqth1}) goes to  $ \int_{0}^{R^{\frac{p-N}{p-1}}} \sqrt[p]{h} \, dt  >0$ as $a \to \infty$, which contradicts $(\ref{sqth1}).$ Thus, $v_a$ must have a first zero, $z_{a},$ with $M_a<z_a<R^{\frac{p-N}{p-1}}$ and $v_a'< 0$ on $(M_a, z_a].$\\
Now, integrating $(\ref{rewrite001})$ on $(M_a, z_{a})$ gives
\begin{equation}
    \int_{M_a}^{z_{a}} \sqrt[p]{h}  \, ds  \leq \sqrt[p]{\frac{p-1}{p}}\int_{0}^{v_a(M_a)}  \frac{  dt}{  \sqrt[p]{F(v_a(M_a)) - F(t)}  } \label{hf001}.
\end{equation}
Again, since $f$ is superlinear, the right-hand side of $(\ref{hf001})$ goes to zero as $a\to \infty$, therefore $z_{a} - M_a \to 0$ as $a \to \infty.$ Since $M_a \to 0$ as $a \to 0$, by Lemma 2.5 it follows that $z_{a} = (z_{a} -M_a) + M_a \to 0$ as $ a \to \infty$.  \qed

\vskip .2 in

{\bf Lemma 2.7:} Let $N>2, a>0$, and assume {\it (H1)}-{\it(H4)}. Suppose $v_a$ solves (\ref{e7})-(\ref{e8}). Let $0< z_a <R^{\frac{p-N}{p-1}}$ and $v_a(z_a)=0$. Then $|v_a'(z_a)|\to \infty$ as $a\to \infty.$

\vskip .1 in

{\bf Proof:} Suppose $v_a(z_a)=0$ for some $0< z_a <R^{\frac{p-N}{p-1}}.$ Then there exists a local extrema on $(0, z_a)$, so suppose without loss of generality there is a local maximum, $M_a,$ with $0<M_a<z_a$ such that $v_a'(M_a)=0,$ and $v_a'(t)< 0$ for all $t\in (M_a, z_a].$ \\
From Lemma 2.5, we know that $v_a(M_a)\to \infty$ as $a\to \infty.$ Therefore, it follows that $F(v_a(M_a))\to \infty$ as $a\to\infty.$ Also evaluating $(\ref{frutti})$ at $z_a$ and recalling that $h$ is decreasing, we obtain
$$ |v_{a}'(z_a)|^p\geq \frac{p}{p-1} h(z_a)F(v_a(M_a))\geq  \frac{p}{p-1} h(R^{\frac{p-N}{p-1}})F(v_a(M_a))\to \infty \quad \text{as} \ a\to \infty.$$
Therefore, $|v_a'(z_a)|\to \infty$ as $a\to \infty.$ \qed

\vskip .3 in

We may now repeat this argument on $(z_a, R^{\frac{p-N}{p-1}})$ and show that $v_a$ has as many zeros as desired as $a \to \infty$ and at each zero, $z_a,$ we have $|v_a'(z_a)|\to \infty$ as $a\to \infty.$

\vskip .2 in

\section{Proof of Theorem 1}
Let $a>0$, and $n$ be a non-negative integer. Now let
$$S_{n} = \{ a > 0  \ | \ v_a \textrm{ has exactly } n \textrm{ zeros  on } (0, R^{\frac{p-N}{p-1}}) \}.$$
From Lemma 2.6, we know some of the $S_{n}$ are non-empty so let $n_0$ be the smallest value of $n$ such that $S_{n_0}\neq \emptyset.$ It follows from Lemma 2.3 that $S_{n_0}$ is bounded from above. Hence let
$$a_{0}:=\sup S_{n_0}.$$
Now we prove that $v_{a_{0}}$ has $n_0$ zeros on $(0, R^{\frac{p-N}{p-1}}).$ For $a<a_{0},$ $v_a$ has exactly $n_0$ zeros on $(0, R^{\frac{p-N}{p-1}})$ and at each zero, $z,$ we have $v_a'(z)\neq 0$ by Lemma 2.2. Therefore, by continuity with respect to $a$ we see that $v_{a_{0}}$ has at least $n_0$ zeros on $(0, R^{\frac{p-N}{p-1}}).$ If $v_{a_{0}}$ has an $(n_0+1)^{st}$ zero on $(0, R^{\frac{p-N}{p-1}}),$ then so would $v_a$ for $a$ slightly smaller than $a_{0}$ (by continuity with respect to $a$) contradicting the definition of $a_{0}.$ Thus $v_{a_{0}}$ has exactly $n_0$ zeros on $(0, R^{\frac{p-N}{p-1}}).$ Denote $z_{n_0}$ as the $n_0^{th}$ zero on $(0, R^{\frac{p-N}{p-1}})$ and without loss of generality let us assume $v_{a_0}>0$ on $(z_{n_0}, R^{\frac{p-N}{p-1}})$. Then by the continuity of $v_{a_0}$ it follows that $v_{a_0}(R^{\frac{p-N}{p-1}})\geq 0.$ Now if $v_{a_0}(R^{\frac{p-N}{p-1}})> 0,$ then by continuous dependence, and since $v_a$ has exactly $n_0$ zeros on $(0, R^{\frac{p-N}{p-1}})$ for $a>a_0$ and $a$ close to $a_0$ , we get $v_{a}(R^{\frac{p-N}{p-1}})> 0$ for $a>a_0$ and $a$ close to $a_0,$ a contradiction. Therefore $v_{a_{0}}(R^{\frac{p-N}{p-1}})=0$, and hence $v_{a_0}$ is a solution of $(\ref{e9})$-$(\ref{e10})$ with exactly $n_0$ zeros. In addition, by Lemma 2.2, $v_{a_{0}}'(R^{\frac{p-N}{p-1}})\neq 0.$ 
\vskip .1 in
Now since $v_{a_{0}}'(R^{\frac{p-N}{p-1}})\neq 0,$ it follows that if $a$ is slightly larger than $a_{0}$ then $v_a$ has exactly $n_0+1$ zeros on $(0, R^{\frac{p-N}{p-1}}).$ Thus, $S_{n_0+1}$ is non-empty. Let
$$a_1:=\sup S_{n_0+1}.$$
In a similar way, we can show $v_{a_1}$ has exactly  $n_0+1$ zeros on $(0, R^{\frac{p-N}{p-1}})$ and $v_{a_1}(R^{\frac{p-N}{p-1}})= 0, \ v_{a_1}'(R^{\frac{p-N}{p-1}})\neq 0.$   
\vskip .1  in
Continuing in this way, we see that for every $n\geq n_0$ there exists a solution of (\ref{e7})-(\ref{e8}) with exactly $n$ zeros on $(0, R^{\frac{p-N}{p-1}})$ and $v_{a_n}(R^{\frac{p-N}{p-1}})= 0.$ This completes the proof. \qed

\vskip .2 in

P.O. Box 311430

Department of Mathematics

University of North Texas

Denton, TX 76203-1430
\vskip .1 in

mdsuzanahamed@my.unt.edu

iaia@unt.edu

\end{document}